\documentclass[oneside]{amsart}

\usepackage[textwidth=13.8cm]{geometry}
\usepackage{amsmath, amssymb, amsthm}
\usepackage{pstricks, pst-node}
\usepackage{cite}
\usepackage{graphicx}
\usepackage{pinlabel}
\usepackage{enumitem}
\usepackage{verbatim}
\usepackage{epigraph}
\usepackage{tikz}
\usetikzlibrary{matrix,arrows,decorations.pathmorphing}
\usepackage{longtable}

\newtheorem{theorem}{Theorem}
\newtheorem{lemma}[theorem]{Lemma}
\newtheorem{prop}[theorem]{Proposition}

\theoremstyle{definition}

\theoremstyle{remark}
\newtheorem*{remark}{Remark}

\numberwithin{equation}{section}

\newcommand{\R}{\mathbb{R}}
\newcommand{\Z}{\mathbb{Z}}

\newcommand{\bdry}[1]{\partial_{\infty}#1}

\newcommand{\h}{\mathfrak{h}}

\newcommand{\modulus}[1]{\left| #1 \right|}

\newcommand{\T}{\mathbb{T}}

\newcommand{\CLF}{\mathrm{CLF}}

\newcommand{\DL}{\mathrm{DL}}

\begin{document}
	
\title{Geometry of the conjugacy problem in lamplighter groups}

\author{Andrew Sale}
\address{Department of Mathematics, 1326 Stevenson Center, Vanderbilty University, Nashville TN 37240, USA}

\email{andrew.sale@some.oxon.org}


\begin{abstract}
	In this note we investigate the conjugacy problem in lamplighter groups with particular interest in the role of their geometry. In particular we show that the conjugacy length function is linear.
\end{abstract}

\maketitle

The conjugacy (search) problem of Max Dehn is a century-old algorithmic question in group theory that has received much attention of late due to potential applications to cryptography, see for example \cite{AAG99,KLCHKP00}.
The conjugacy problem asks whether, given a group $\Gamma$ with finite generating set $A$, there exists an algorithm which, on input words $u,v$ on $A\cup A^{-1}$, determines whether $u,v$ represent conjugate elements in $\Gamma$.
The conjugacy search problem is similarly themed, but instead the input is two elements known to be conjugate, and the algorithm should produce a conjugating element.

We study the \textit{conjugacy length function}, a quantitative version of the conjugacy problem. 
Given a group $\Gamma$ and a length function $\left| \cdot \right|:\Gamma\to [0,\infty),$ (typically word length) the conjugacy length function
is the minimal function $$\CLF_\Gamma : [0,\infty) \rightarrow [0,\infty)$$ which satisfies the following: if $u$ is conjugate to $v$ in $\Gamma$ and $\modulus{u}+\modulus{v}\leq n$ then there exists a conjugator $\gamma \in \Gamma$ such that $\modulus{\gamma} \leq \CLF_\Gamma(n)$. This function has been estimated for many classes of groups, a brief list includes hyperbolic groups \cite{lysenok} and mapping class groups \cite{MM00,BD11,Tao11} where in both it is known to be (at most) linear, it is bounded by an exponential function for CAT(0) groups \cite{BH99}, and in free solvable groups it is at most cubic \cite{sale2013wreath}. 

In \cite{sale2013wreath} the author investigated the behaviour of the conjugacy length function under wreath products. For the class of lamplighter groups $\Z_q \wr \Z$ or $\Z\wr\Z$, this result is not optimal, giving a cubic upper bound. We show here that it is linear. In fact, we consider the Diestel--Leader groups $\Gamma_2(R)$, which are defined in Section \ref{sec:DL groups and gen lamplighters}, that includes the lamplighter groups $\Gamma_2(\Z_q)=\Z_q\wr\Z$ and $\Gamma_2(\Z)=\Z\wr\Z$.

\begin{theorem}\label{thm:clf for lamplighter}
	Let $R$ be a commutative ring with unity and let $\Gamma = \Gamma_2(R)$. Then there is a generating set $S$, which is finite when $R$ is finite, such that with respect to the corresponding word length we have
	$$\CLF_\Gamma(n) \leq 3n.$$
\end{theorem}

A key feature we use is the fact that the Cayley graphs of the groups $\DL_2(R)$ are the horocyclic product of two trees \cite{BNW08,AR14}, a special case of a Diestel--Leader graph, introduced in \cite{DL01}. The motivation for the geometric approach that we use is so that the method may generalise to a wider class of groups, the higher-rank Diestel--Leader groups $\Gamma_d(R)$ for $d \geq2$. Indeed, C.{} Abbott has done precisely this, obtaining an exponential upper-bound for their conjugacy length functions  when $R$ is finite \cite{Abbott}.

The tools developed for Theorem \ref{thm:clf for lamplighter} are then applied to give an algorithm solving the conjugacy search problem. The algorithm will run in quadratic time with respect to the length of the input words when $R$ is finitely generated as an abelian group.

Theorem \ref{thm:clf for lamplighter} originally appeared with an algebraic proof in a preprint of the author \cite{Sale11}, which was part of the inspiration for \cite{Sale12groupext}. The latter looks at the behaviour of conjugacy length under group extensions, and the proof of Theorem \ref{thm:clf for lamplighter} could be rewritten in the language used therein. This note has been modified from the geometric proof in the author's doctoral thesis \cite{SaleThesis}.

The structure is as follows. We begin by defining the groups $\Gamma_d(R)$ in Section \ref{sec:lamplighter groups} and explain their geometry through Diestel--Leader graphs.
A word length estimate, which follows from a formula of Stein and Taback \cite{ST12}, is given in Section \ref{sec:word length}, while Theorem \ref{thm:clf for lamplighter} is proved in Section \ref{sec:clf} and the algorithm is described in Section \ref{sec:algorithm}.

\section{Lamplighter groups as Diestel--Leader groups}\label{sec:lamplighter groups}

\subsection{Horocyclic products, Diestel--Leader graphs and $R$--branching trees}\label{sec:Lamplighers:horocyclic products}

We give here a brief introduction to horocyclic products and Diestel-Leader graphs. For a more complete description see for example \cite{BNW08}. 

Let $T$ be a simplicial tree and $\omega \in \bdry{T}$ a boundary point of $T$, that is $\omega$ is an equivalence class of asymptotic geodesic rays.
Recall that in $T$, two geodesic rays $\rho_i:[0,\infty) \to T$, for $i=1,2$, are asymptotic if and only if they merge: there exists $x>0$ and $s\in \R$ such that $\rho_1(t)=\rho_2(t+s)$, for all $t\geq x$.

For any vertex $x$ in $T$ there is a unique geodesic ray emerging from $x$ that is in $\omega$. 
Given a pair of vertices, $x,y$, the corresponding rays will coincide from some vertex $x \curlywedge y$ onwards. We call this the \emph{greatest common ancestor} of $x$ and $y$.

After fixing a basepoint $o$ in the vertex set of $T$ we can define a Busemann function $\mathfrak{h}$ on the vertices of $T$ as
$$
\mathfrak{h}(y) = d_T(y,o \curlywedge y) - d_T(o , o \curlywedge y) \textrm{.}$$
The difference in the value of the Busemann function gives the distance from a vertex $x$ to its common ancestor with some other vertex $y$. That is
		$$d_T(x,x\curlywedge y) = \h(x)-\h(x \curlywedge y).$$
The $k$--th horocycle of $T$ based at $\omega$ is $H_k = \{y \in T \mid \mathfrak{h}(y)=k\}$.

Given a collection $T_1 , \ldots , T_n$ of simplicial trees together with a chosen collection of respective Busemann functions $\mathfrak{h}_1 , \ldots , \mathfrak{h}_n$, we define the \emph{horocyclic product} to be
\begin{equation}\label{eq:horocyclic product}
\prod_{i=1}^n {}_\mathfrak{h} T_i = \left\{ (y_1, \ldots , y_n) \in T_1 \times \ldots \times T_n \mid \sum_{i=1}^n \mathfrak{h}_i(y_i) = 0 \right\}\textrm{.}
\end{equation}

The \emph{Diestel--Leader graph} $\mathrm{DL}(q_1 , \ldots , q_d)$ is the $1$--skeleton of the horocyclic product of trees $\mathbb{T}_{q_1},\ldots,\T_{q_d}$, where $\mathbb{T}_q$ is the $q+1$ regular tree. When $q_1=q_2=\ldots =q_d=q$, the corresponding Diestel-Leader graph is also denoted by $\mathrm{DL}_d(q)$.

Following \cite{AR14}, we can extend this to the notion of an $R$--branching tree, where $R$ is any commutative ring with unity. Such a tree, denoted $\T_R$, has the property that for every vertex $v$, the set $E(v)$ of edges touching $v$ is in a bijection with $\{x_0\}\cup R$, where $x_0\notin R$. We denote by $\DL_d(R)$ the horocyclic product of $d$ copies of $\T_R$.

\subsection{Diestel--Leader groups and generalised lamplighters}\label{sec:DL groups and gen lamplighters}

The fact that the Diestel--Leader graph $\DL_2(q)$ is a Cayley graph for the lamplighter group $\Z_q \wr \Z$ was explained in \cite{Woes05}. This is a special case of the following result of Bartholdi, Neuhauser and Woess:

\begin{theorem}[Bartholdi--Neuhauser--Woess {\cite[(3.14)]{BNW08}}]\label{theorem:BNW lamplighter cayley graph}
The Diestel--Leader graph $\DL_d(q)$ is a Cayley graph of a group, denoted $\Gamma_d(R)$, where $R$ is a commutative ring with unity of order $q$.
\end{theorem}

The groups $\Gamma_d(R)$ are called \emph{Diestel--Leader groups}. When $d=2$ they include the lamplighter groups $\Z_q\wr \Z$, and when $d=3$ we obtain groups previously considered by Baumslag \cite{Baum72}, \cite{Baum74}. The descriptions given below for the lamplighter groups can be extended to the Diestel--Leader groups of more than $2$ trees. The word length in Diestel--Leader groups is studied by Stein and Taback \cite{ST12} where they give a formula for its calculation.

As described by Amchislavska and Riley \cite[Section 1.4]{AR14} we can consider groups whose Cayley graph is the horocyclic product of $R$--branching trees $\T_R$. Mixing the notation of \cite{AR14} and \cite{BNW08}, we let $\Gamma_d(R)$ denote the group whose Cayley graph is the horocyclic product of $d$ trees $\T_R$.
In this note will focus on the case when $d=2$. We can recognise these as groups of the form
		\begin{equation}
		\label{eqb:semidirect}
		\Gamma_2(R) \cong R[t,t^{-1}]\rtimes \Z,
		\end{equation} 
where $\Z$ acts on $R[t,t^{-1}]$ by multiplication by $t$.
We refer the reader to \cite{AR14} for more details, particularly for when $d>2$, where these groups are denoted by $\Gamma_{d-1}(R)$ --- their subscript agrees with the number of copies of $\Z$ in the quotient.

The group $\Gamma_2(R)$  can be represented by matrices
\begin{equation*}
\left\{ \left(\begin{array}{cc}t^k & P \\ 0 & 1\end{array}\right) : k \in \Z \textrm{, } P \in R [t^{-1},t] \right\}\mathrm{,}
\end{equation*}
with generating set
\begin{equation*}
S_0=\left\{ \left(\begin{array}{cc}t & b \\ 0 & 1\end{array}\right) : b \in R \right\}\mathrm{.}
\end{equation*}
For the sake of notation, we use the identification from \eqref{eqb:semidirect}:
		$$\left(\begin{array}{cc}t^k & P \\ 0 & 1\end{array}\right) \mapsto (P,k).$$	
Amchislavska and Riley generalised Theorem \ref{theorem:BNW lamplighter cayley graph}, giving

\begin{theorem}[Amchislavska and Riley \cite{AR14}]\label{thm:margarita-tim}
	Let $R$ be a ring with unity.
	There is a generating set $S$ for $\Gamma_d(R)$, with respect to which the Cayley graph is $\DL_d(R)$.
	
	When $d=2$, the generating set is $S=S_0$.
\end{theorem}
When $R$ is infinite, the generating set $S$ is infinite too. This is to be expected, since the valency of the trees, and hence the Diestel--Leader graph, are infinite.

\subsection{Realising $\DL_2(R)$ as the Cayley graph of a lamplighter group}\label{sec:DL to lamplighter}

We begin by labelling the edges of $\T_R$ with elements from $R$.
We pick a bi-infinite geodesic $\rho:\R \to \T_R$ so that its restriction to $[0,\infty)$ is in $\omega\in\bdry{\T_R}$.
Let $o=\rho(0)$ and consider the corresponding Busemann function $\h$.
Recall that at each vertex $v$ we have a bijection 
		\begin{equation}\label{eqn:vertex bijection}
		E(v) \to \{x_0\} \cup R
		\end{equation}
where $E(v)$ is the set of edges touching $v$, and $x_0\notin R$.
At each vertex $v$ there is a unique edge $e_v$ in $E(v)$ which forms the first edge in a geodesic ray asymptotic to $\rho$.
The labelling is an iterative method, as follows.
First label each edge in $\rho$ by $0$.
Then for vertices $v$ in which $e_v$ has been previously labelled, but some edges in $E(v)$ remain unlabelled, use the bijection \eqref{eqn:vertex bijection} to label those remaining edges.
Repeat this step.

We now take two copies of our decorated tree, $T_1,T_2$, with boundary points $\omega_1,\omega_2$ and geodesics $\rho_1,\rho_2$ respectively, and take their horocyclic product to get $\DL_2(R)$. 
Given a vertex $x=(x_1,x_2)$ of $\DL_2(R)$ we consider the two geodesic rays starting at $x_1,x_2$ respectively that are in $\omega$. 
Similar to \cite[Lemma 3.1]{AR14} we get a finitely supported bi-infinite sequence $(a_i)$ in $R$, obtained by reading the edge labels in these two rays as follows.
For each $i$ there is exactly one edge in the rays that travels between vertices in $H_{i-1}$ and $H_{i}$ if in $T_1$, and $H_{-i+1}$ and $H_{-i}$ if in $T_2$. The label of this edge determines $a_i$.
Thus, $(a_{\h(x_1)},a_{\h(x_1)-1},\ldots)$ is the sequence of edge labels of the ray in $T_1$, and $(a_{\h(x_1)+1},a_{\h(x_1)+2},\ldots)$ is the sequence from the ray in $T_2$.

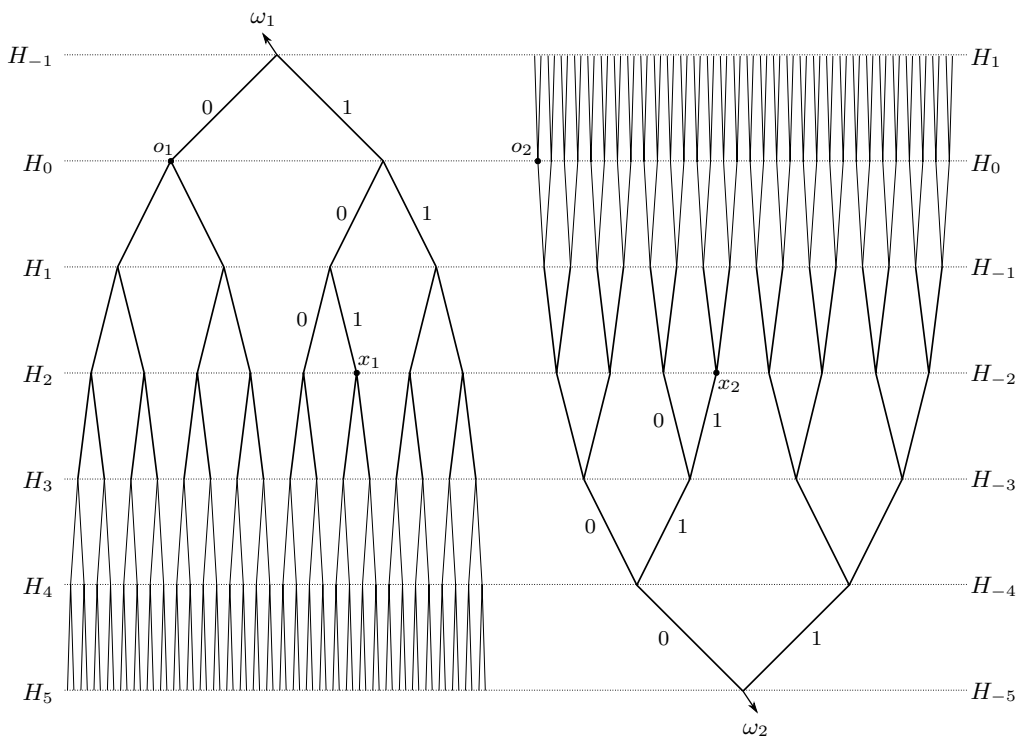
\begin{figure}[h!]
\labellist \small\hair 4pt
		\pinlabel $\omega_1$ [b] at 120 409
		\pinlabel $\omega_2$ [t] at 417 0
		\pinlabel $H_{-1}$ [r] at 0 396
		\pinlabel $H_0$ [r] at 0 332
		\pinlabel $H_1$ [r] at 0 268
		\pinlabel $H_2$ [r] at 0 204
		\pinlabel $H_3$ [r] at 0 140
		\pinlabel $H_4$ [r] at 0 76
		\pinlabel $H_5$ [r] at 0 12
				\pinlabel $H_{1}$ [l] at 540 396
				\pinlabel $H_0$ [l] at 540 332
				\pinlabel $H_{-1}$ [l] at 540 268
				\pinlabel $H_{-2}$ [l] at 540 204
				\pinlabel $H_{-3}$ [l] at 540 140
				\pinlabel $H_{-4}$ [l] at 540 76
				\pinlabel $H_{-5}$ [l] at 540 12
	\footnotesize
		\pinlabel $0$ [r] at 96 366
		\pinlabel $1$ [l] at 161 366
			\pinlabel $0$ [r] at 177 302
			\pinlabel $1$ [l] at 209 302
		\pinlabel $0$ [r] at 153 238
		\pinlabel $1$ [l] at 167 238
			\pinlabel $0$ [r] at 368 177
			\pinlabel $1$ [l] at 384 177
		\pinlabel $0$ [r] at 327 113
		\pinlabel $1$ [l] at 363 113
			\pinlabel $0$ [r] at 371 45
			\pinlabel $1$ [l] at 444 45
		\footnotesize\hair 4pt
			\pinlabel $x_1$ [b] at 184 202
			\pinlabel $x_2$ at 401 198
			\pinlabel $o_1$ at 60 341
			\pinlabel $o_2$ at 276 341
\endlabellist
 \centering\hspace{0.5cm}\includegraphics[width=12cm]{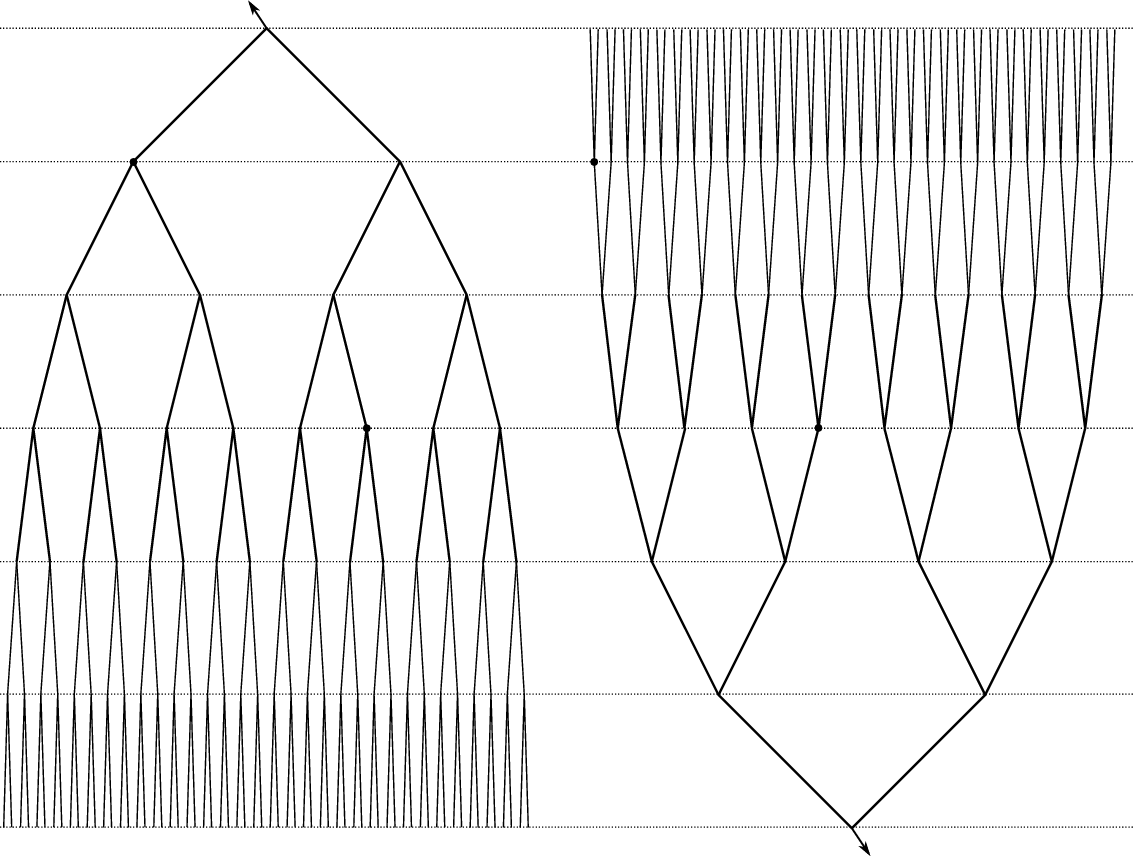}
\caption{The tree $T_1$ is on the left, and $T_2$ on the right. 
	From $T_1$ we get $a_2=1,a_1=0,a_{0}=1$ and $a_i=0$ for $i<0$. From $T_2$ we get $a_3=1,a_4=1$ and $a_i=0$ for $i>4$. Thus $f_x = 1+t^2+t^3+t^4$. }\label{fig:DL to poly}
\end{figure}

We can use the sequence $(a_i)$ to determine an element of $R[t,t^{-1}]$:
	\begin{equation}\label{defn:f_x}
	f_x:=\sum\limits_{i=-\infty}^\infty a_it^i.
	\end{equation}
An example is given in Figure \ref{fig:DL to poly}. We use this to establish an identification
		\begin{equation} \label{eq:Cayley graph identification}
		\DL_2(R) \ni x = (x_1,x_2) \longleftrightarrow (f_x, \h(x_1))\in \Gamma_2(R).
		\end{equation}
Note that comparing this identification with that in \cite[Section 4]{AR14}, the role of our trees $T_1,T_2$ are exchanged.

Let $f \in R[t,t^{-1}]$.
There exists some $n\in\Z$ such that the coefficient of $t^k$ is zero for all $k < n$. The valuation $v_0(f)$ is defined to be the supremum of all such $n$ (note that we define $v_0(0)=\infty$).
Similarly we set $v_0^-(f)$ to be the largest $n$ for which the coefficient of $t^n$ in $f$ is non-zero.
Geometrically, $v_0(f)-1$ gives the horocycle in which the geodesic ray emerging from $x_1$, asymptotic to $\omega_1$, merges with $\rho_1$. Similarly, $-v_0^-(f)$ give the horocycle in $T_2$  where the geodesic ray in $\omega_2$ that emerges from $x_2$ merges with $\rho_2$.

Before proceeding further, we make the following observation which will be useful when discussing word length.

\begin{lemma}\label{lem:height of ancestor under action}
	Let $x=(x_1,x_2),y=(y_1,y_2)\in \DL_2(R)$ and $\gamma \in \Gamma_2(R)$. Then
	$$\mathfrak{h}(\gamma x_i \curlywedge \gamma y_i)=\mathfrak{h}(x_i \curlywedge y_i)+\h(\gamma o_i).$$
\end{lemma}

\begin{proof}
		We prove this for the first tree. For the second tree the proof is analogous.
		
		Let $f_x=\sum a_it^i$ and $f_y=\sum b_it^i$ be as in \eqref{defn:f_x}. 
		We claim that
			$$\mathfrak{h}(x_1\curlywedge y_1) = \min \{ v_0(f_x-f_y)-1, \mathfrak{h}(x_1),\mathfrak{h}(y_1)\}.$$
		Since both have finite support, there is some $k$ such that for $i<k$ we have $a_i=b_i$. This $k$ is the valuation $v_0(f_x-f_y)$. 
		First suppose that $x_1\curlywedge y_1$ is distinct from $x_1$ and $y_1$.
		Then $v_0(f_x-f_y)$ gives the last horocycle before the merging of the two geodesics in $\omega_1$ emerging from $x_1$ and $y_1$ respectively.
		Hence $\h(x_1\curlywedge y_1) = v_0(f_x-f_y)-1$.
		
		On the other hand, if the common ancestor is one of the two given vertices, say $x_1$, then this valuation may not give the correct horocycle. Indeed it may be that $f_x$ and $f_y$ agree for coefficients of $f_x$ that come from $T_2$ and have no consequence for the location of $x_1$. Hence, in this case we see that we have $\h(x_1\curlywedge y_1) = \h(x_1)$. This proves the claim.

		Suppose $\gamma=(P,s)$. Then $f_{\gamma x} = P+t^sf_x$ and $f_{\gamma y} = P+t^sf_y$.
		We thus obtain, using the above claim, the following for $\mathfrak{h}(\gamma x_1\curlywedge \gamma y_1)$, noting that the terms in $P$ will both cancel, leaving
			$$\mathfrak{h}(\gamma x_1 \curlywedge \gamma y_1) = \min\{v_0(t^sf_x-t^sf_y)-1,\mathfrak{h}(\gamma x_1),\mathfrak{h}(\gamma y_1)\}.$$
		Since $v_0(t^sf_x-t^sf_y)=v_0(f_x-f_y)+s$, $\mathfrak{h}(\gamma x_1)=\mathfrak{h}(x_1)+s$, and similarly for $y_1$, the Lemma holds.
%
%
%
%
%
\end{proof}

We remark that Lemma \ref{lem:height of ancestor under action} can be generalised to $\DL_d(q)$ for $d\geq 2$.

\subsection{The (right) action of generators on $\DL_2(R)$}

	Recall we have two trees $T_1,T_2$, each isomorphic to $\T_R$ and coming with an equivalence class of rays $\omega_i$ in $\bdry{T_i}$.
	Starting from a vertex $x_i \in T_i$, we will use the phrase \textit{``go up in $T_i$''} to mean move along the (unique) edge in $T_i$ which is in the geodesic ray emanating from $x_i$ that is in $\omega_i$. To \textit{``go down in $T_i$''} means we move along any other edge. Recall this edge will be labelled by an element of $R$, as described in Section \ref{sec:DL to lamplighter}.
	
	\begin{lemma}\label{lem:generator to instruction}
		A letter $(b,1)\in S$ in a word is an instruction:
		move up one edge in tree $T_2$, suppose this has label $a$, and move down an edge labelled $a+b\in R$ in tree $T_1$.
	\end{lemma}
	
	\begin{proof}
		We have an edge in the Cayley graph of $\Gamma_2(R)$ labelled by $(b,1)$ starting at $(f,k)$ and finishing at $(f,k)(b,1) = (f+bt^k,k+1)$. 
		Thus, in $T_1$ we travel from a vertex in $H_k$ down to $H_{k+1}$.
		Suppose $a$ is the coefficient of $t^k$ in $f$, then the edge we followed, by the identification \eqref{eq:Cayley graph identification}, has label $a+b$, which is the coefficient of $t^k$ in $f+bt^k$.
		
		Meanwhile, in $T_2$ we have travelled up an edge, from $H_{-k}$ to $H_{-k-1}$. The edge we followed was labelled by the coefficient of $t^k$ in $f$, namely $a$.
	\end{proof}
	
	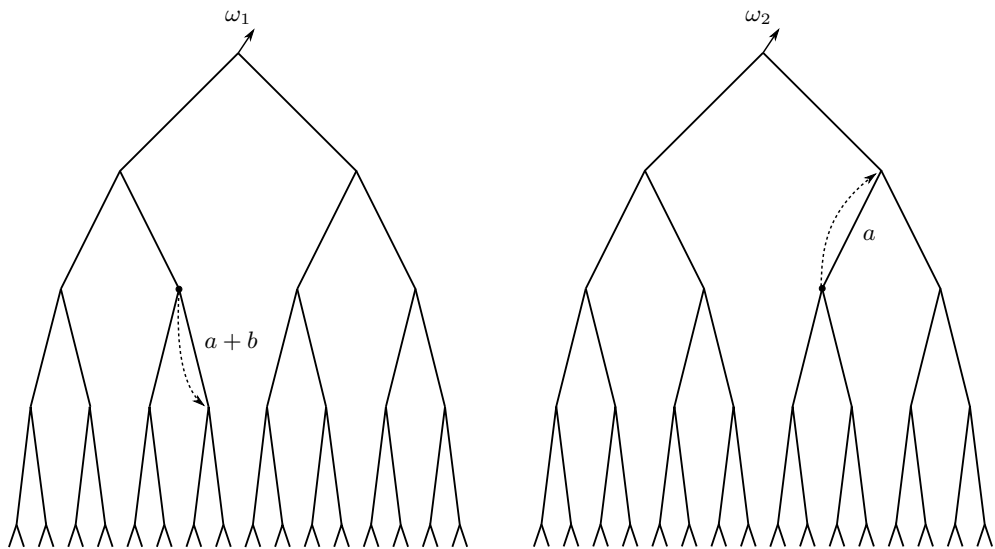
\begin{figure}
		\labellist
		\small\hair 5pt
		\pinlabel $\omega_1$ at 125 288
		\pinlabel $\omega_2$ at 407 288
		\pinlabel $a$ at 467 170
		\pinlabel $a+b$ at 121 112
		\endlabellist
		
		\includegraphics[width=0.95\linewidth]{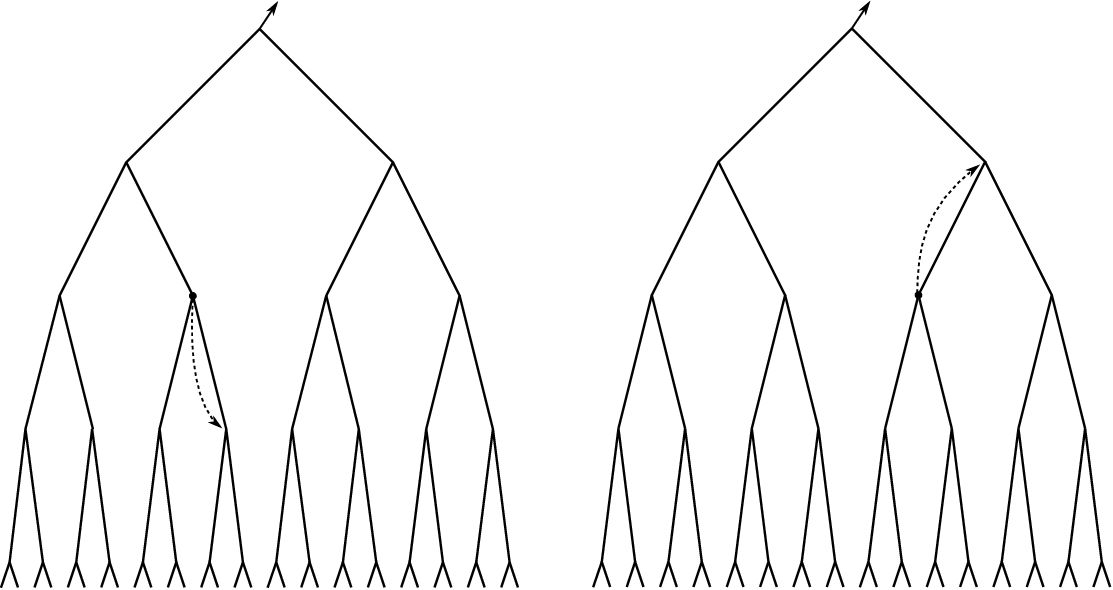}
		\caption{The action of a letter $(b,1)$ within a word tells us to move down in $T_1$ and up in $T_2$. The label of the edge we move up in $T_2$ combines with $b$ to tell us which edge to go down in $T_1$.}
	\end{figure}

\section{Word length}\label{sec:word length}

The word length of elements in a general Diestel--Leader group $\Gamma_d(R)$, when $R$ is finite and $d\geq 2$, has been studied by Stein and Taback \cite{ST12}. They give a formula for the word length of an element by looking at the climb and fall of a geodesic path in each tree from the basepoint to its image under the action of the element. Let $o_i$ denote the basepoint of the $i$--th tree. For $g \in \Gamma_d(R)$, denote by $m_i(g)$ the length of the climb of the geodesic from $o_i$ to $go_i$ and $l_i(g)$ the length of the fall. More concretely we mean:
		$$m_i(g) = d(o_i , o_i \curlywedge go_i), \ \ l_i(g) = d(go_i,o_i \curlywedge go_i).$$
Stein and Taback's formula, see \cite[Section 3]{ST12}, is precise, but can be used to give an estimate for the word length of $g$ that is practical for our purpose.
We include a proof of this estimate not only for completeness, but also because we wish to take $R$ to be finite or infinite.

\begin{prop}[Stein--Taback word length estimate \cite{ST12}]\label{prop:stein-taback}
	Let $g \in \Gamma_2(R)$. Then
		$$m_1(g)+m_2(g) \leq \modulus{g} \leq 2(m_1(g)+m_2(g))$$
where $\modulus{g}$ is the word length of $g$ with respect to the generating set $S$ of Theorem \ref{thm:margarita-tim}.
\end{prop}

\begin{remark}
	We may replace $m_i$ with $l_i$ in the above inequality, since whenever we are in a horocyclic product we will have
			$$
			\sum_{i=1}^d m_i(g) = \sum_{i=1}^d l_i(g).
			$$
\end{remark}

\begin{proof}[Proof of Proposition \ref{prop:stein-taback}]
	With Lemma \ref{lem:generator to instruction} in mind, to reach $(go_1,go_2)$ from $(o_1,o_2)$ in $\DL_2(R)$, we must climb at least $m_1(g)$ in $T_1$ and $m_2(g)$ in $T_2$, thus implying the lower bound.
	
	Without loss of generality we may assume $\h(go_1)\geq 0$.
	First apply ${m_1(g)}$ generators to get to $o_1 \curlywedge go_1$ in $T_1$. 
	Second apply the appropriate sequence of downward movements in $T_1$ to reach $go_1$. 
	This sequence will be a word of length $l_1(g)$. 
	Meanwhile, in $T_2$ we have reached a vertex in the same horocycle as $go_2$. 
	The third step is to climb to $o_2\curlywedge go_2$ and drop down to $go_2$.
	During this third phase, we will drop down in $T_1$, but then climb, inevitably retracing our steps, back up to $go_1$.	
	The third and final step requires $2l_2(g)$ elements from $S$.
	So in total, we have a word of length
			$$m_1(g) + l_1(g) + 2l_2(g).$$
	Since we are in a horocyclic product we have $m_1(g) + m_2(g) = l_1(g) + l_2(g)$. We also have $l_2(g)\leq m_2(g)$ since $\h(go_1)\geq 0$. Hence we get
			$$m_1(g) + l_1(g) + 2l_2(g) \leq 2(m_1(g) +m_2(g))$$
	as required.	
\end{proof}

The notion of climbing and falling can be extended to paths describing the concatenation of words. In particular, for $g,h \in \Gamma_d(R)$ let
$$m_{h,i}(g)= d(ho_i,ho_i \curlywedge go_i) , \ \ l_{h,i}(g)=d(go_i , ho_i \curlywedge go_i).$$
A consequence of Lemma \ref{lem:height of ancestor under action} is that we can measure the length of the word by looking at its action on different points in each tree:

\begin{lemma}\label{lem:m_i=m_h,i}
For every $g,h \in \Gamma_2(R)$ and each $i=1,2$ we have the following:
		$$m_i(g)=m_{h,i}(hg), \  \ \ l_i(g)=l_{h,i}(hg).$$
\end{lemma}

\begin{proof}
By Lemma \ref{lem:height of ancestor under action}, in each tree, the height of the common ancestor $o_i \curlywedge g o_i$ never travels too far from the height of $ho_i \curlywedge hg o_i$. To be precise:
		$$\mathfrak{h}_i(o_i \curlywedge g o_i) = \mathfrak{h}_i(ho_i \curlywedge hg o_i) - \mathfrak{h}_i(ho_i)$$
Hence 
\begin{eqnarray*}
m_i(g) & = & d(o_i,o_i \curlywedge go_i)\\
	& = & -\mathfrak{h}_i(o_i \curlywedge go_i)\\
	& = & -\mathfrak{h}_i(ho_i \curlywedge hgo_i)+ \mathfrak{h}_i(ho_i).
\end{eqnarray*}
But $d(ho_i,ho_i \curlywedge hgo_i) = \mathfrak{h}_i(ho_i) - \mathfrak{h}_i(ho_i\curlywedge hgo_i)$, so we get $m_i(g)=m_{h,i}(hg)$. The result for $l_i(g)$ follows from the result for $m_i(g)$, the relationships:
		$$m_i(g)-l_i(g)=\mathfrak{h}_i(go_i) \ \ \mathrm{and} \ \ m_{h,i}(hg)-l_{h,i}(hg)=\mathfrak{h}_i(hgo_i)-\mathfrak{h}_i(ho_i)$$
and the fact that $\mathfrak{h}_i(go_i)=\mathfrak{h}_i(hgo_i)-\mathfrak{h}_i(ho_i)$.
\end{proof}

Using Lemma \ref{lem:m_i=m_h,i} we can deduce the following estimate from Proposition \ref{prop:stein-taback}:
\begin{equation}\label{eq:lamplighter word length}
\sum\limits_{i=1}^d m_{h,i}(g) \leq \modulus{g^{-1}h} \leq 2\sum\limits_{i=1}^d m_{h,i}(g).
\end{equation}
In these bounds for $\modulus{g^{-1}h}$ we may replace $m_{h,i}$ with $l_{h,i}$ since their sums are equal.
To prove the estimate \eqref{eq:lamplighter word length}, we need
	$$l_i(g^{-1}h) = l_{g,i}(h) = d(ho_i,go_i\curlywedge ho_i) =  m_{h,i}(g)$$
where the first equality is an application of Lemma \ref{lem:m_i=m_h,i}.

\section{The conjugacy length function}\label{sec:clf}

We now prove that the conjugacy length function of $\Gamma_2(R)$ is linear with respect to the word length given by the generating set $S$ defined in Section \ref{sec:DL groups and gen lamplighters}.
The following tells us the structure of short conjugators, from which we deduce \mbox{Theorem \ref{thm:clf for lamplighter}.}

\begin{prop}\label{prop:short conjugators}
	Suppose $u=(P,s),  v=(Q,r) \in \DL_2(R)$. If $u$ is conjugate to $v$ then $r=s$ and either
		\begin{enumerate}
			\item \label{part:unipotent} if $r=s=0$ then there is a conjugator $\gamma=(0,k)$ with $\modulus{k}\leq \modulus{u}+\modulus{v}$;
			\item  if $r=s\neq 0$ then there is a conjugator $\gamma=(f,k)$ with $0\leq k < \modulus{r}$ and
					\begin{equation}\label{eqn:climbing limit}
					m_i(\gamma) \leq \max \{ m_i(u),m_i(v)\}+\modulus{r}.
					\end{equation}
		\end{enumerate}
\end{prop}

\begin{proof}
	Suppose $\gamma=(f,k)$ is a conjugator. By direct calculation, $\gamma  u=v \gamma $ if and only if the following equations hold:
	\begin{eqnarray}
	\label{eq:lamplighter Z}  s+k &=& k+r, \\
	\label{eq:lamplighter functions}  P +t^s f &=& f + t^k Q.
	\end{eqnarray}
From the first of these equations we get $r=s$, and we thus split into the two cases.

\medskip

\textbf{{Case 1}}: $r=0$.
This corresponds to the case when $u$ and $v$ map each basepoint to a vertex of the same height. Equation (\ref{eq:lamplighter functions}) becomes $P=t^kQ$, so we may set $f=0$. This means that $\gamma$ will act on each tree by either a sequence of consecutive up movements or a sequence of consecutive down movements, but never a mixture of both.

As long as $v$ is non-trivial, we may assume that $vo_i \neq o_i$ 
for some $i$.
Assume the action of $\gamma$ on the $i$--th tree corresponds to a consecutive list of downward movements. This case is depicted in Figure \ref{fig:unipotent climb fall} (a).
Then $v\gamma o_i$ will be sitting below $vo_i$, to be precise we have $vo_i \curlywedge v\gamma o_i = vo_i$.
Then the path from $\gamma o_i$ to $\gamma uo_i$ must pass through $vo_i$. Hence $m_{\gamma,i}(\gamma u) \geq \modulus{k}$. Using equation \eqref{eq:lamplighter word length} we then obtain
	\begin{equation*}
		\modulus{u}
		=\modulus{(\gamma u)^{-1}\gamma} 
		\geq m_{\gamma,i}(\gamma u) 
		\geq \modulus{k}.
	\end{equation*}

\begin{figure}[h!]
	\vspace{8mm}
	\labellist \small\hair 5pt
	\pinlabel (a) at 150 360
	\pinlabel (b) at 565 360
	\pinlabel $\gamma o_i$ [r] at 2 2
	\pinlabel $o_i$ [r] at 122 242
	\pinlabel $u$ at 160 211
	\pinlabel $vo_i$ [l] at 203 242
	\pinlabel $\gamma uo_i$ [l] at 227 2
	\pinlabel $o_i$ at 400 2
	\pinlabel $\gamma o_i$ at 515 242
	\pinlabel $\gamma uo_i$ at  650 242
	\pinlabel $vo_i$ at 665 2
	\pinlabel $v$ at 577 211
	\endlabellist
		
	\centering\includegraphics[width=3.2cm,height=40mm]{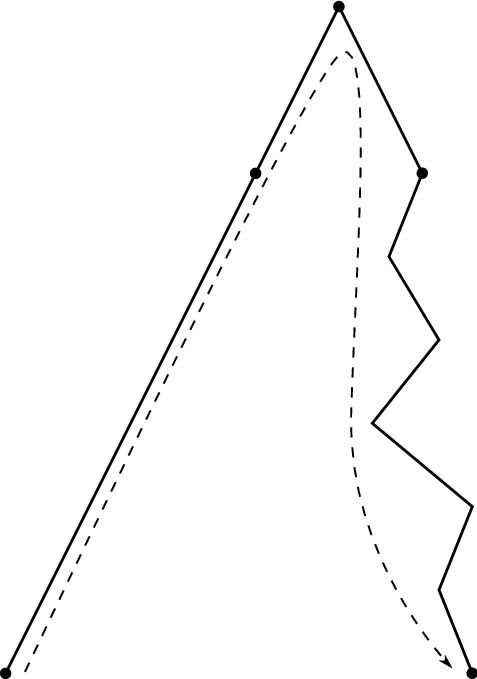}\hspace{25mm}
	\centering\includegraphics[width=3.2cm,height=40mm]{unipotent_climb_fall}
	\caption[Case 1 of Theorem \ref{theorem:clf for lamplighter}]{(a) The action of $\gamma$ is purely downward; (b) the action of $\gamma$ is purely upward.}\label{fig:unipotent climb fall}
\end{figure}

If $\gamma$ acts on the $i$--th tree in a purely upwards manner, then we will get a similar picture (see Figure \ref{fig:unipotent climb fall} (b)), but $\gamma u o_i$ will be sitting directly above $vo_i$, and, provided $\gamma o_i \neq \gamma u o_i$, any path from $o_i$ to $vo_i$ must pass through $\gamma u o_i$.
Here we will get 
		$$\modulus{k} \leq m_i(v) \leq \modulus{v}.$$

The situation that remains is when, up to swapping the trees round, we have $\gamma o_1=\gamma uo_1$, and $o_2=vo_2$, and $\gamma$ acts in a purely upward direction on $T_1$ and downward on $T_2$.
When moving from $v(o_1,o_2)$ to $v\gamma(o_1,o_2)$ in $\DL_2(R)$ we get a pair of geodesics, one in each tree. 
The geodesic in $T_1$ will be from $vo_1$ to $v \gamma o_1 = \gamma u o_1$ and will merge with $\rho_1$ at $o_1 \curlywedge vo_1$, which is in the horocycle at level 
		$$\h(o_1 \curlywedge vo_1) = -m_1(v).$$
Meanwhile, the geodesic in $T_2$ will separate from $\rho_2$, since $\gamma uo_2\neq \gamma o_2$ as otherwise $u$ would be the identity. 
Furthermore, this separation must occur before the geodesic in $T_1$ joins $\rho_1$, since once it has joined $\rho_1$, Lemma \ref{lem:generator to instruction} tells us that in $T_2$ we must follow edges labelled $0$, which would cause us to remain on $\rho_2$ and end up with $\gamma o_2 = \gamma u o_2$.
In $T_2$ the geodesic separates from $\rho_2$ in the horocycle of level $\modulus{k}-m_2(u)$. Hence, by the above argument we must have
		$$\modulus{k} - m_2(u) < m_1(v) \implies \modulus{k} < m_1(u) + m_2(v)\leq \modulus{u}+\modulus{v}.$$
		
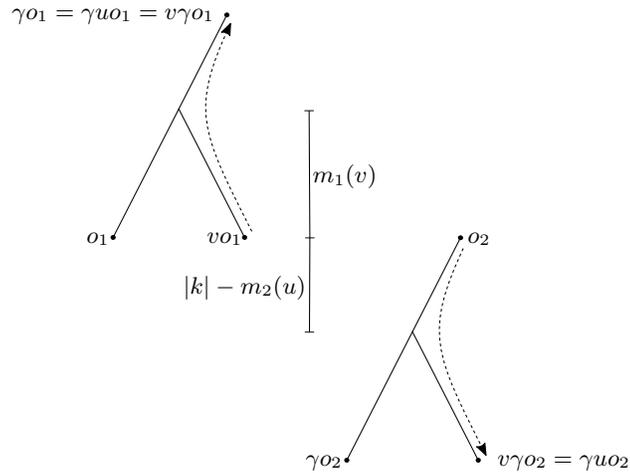
\begin{figure}[h!]
	\labellist \small \hair 2pt
	\pinlabel $o_1$ at -10 185
	\pinlabel $vo_1$ at 88 185
	\pinlabel $\gamma o_1=\gamma uo_1=v\gamma o_1$ [r] at 84 370
	\pinlabel $m_1(v)$ at 181 237
	\pinlabel $o_2$ at 285 185
	\pinlabel $\gamma o_2$ at 165 0
	\pinlabel $v\gamma o_2=\gamma uo_2$ [l] at 295 0
	\pinlabel $\modulus{k}-m_2(u)$ [r] at 155 145
	\endlabellist
	\includegraphics[width=5cm,height=6cm]{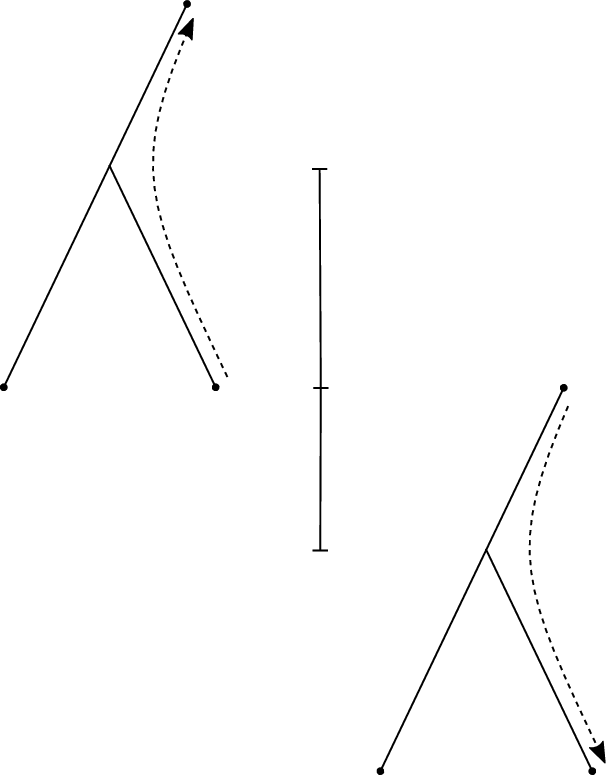}
\caption{The path in $T_1$, on the left, merges with $\rho$ after the path in $T_2$ separates from $\rho$.}\label{fig:special case}
\end{figure}

\medskip
\textbf{{Case 2}}: $r \neq 0$.
By exchanging $u$ and $v$ with their inverses if necessary, we may assume that $r>0$. 
The important step here is to pick the right conjugator. 
Take any conjugator $\gamma'$, satisfying $\gamma' u=v \gamma'$, and premultiply it by a suitable power of $u$ so that we obtain an element $\gamma=u^m\gamma'$, written as above, with $0 \leq k < r$. 
To prove the proposition, we just need to demonstrate the bound \eqref{eqn:climbing limit} on $m_i(\gamma)$. We will show that if this bound were not true then $\gamma u o_i$ and $v\gamma o_i$ would have to be on different branches of their respective trees, see Figure \ref{fig:CLF lamp both trees}.

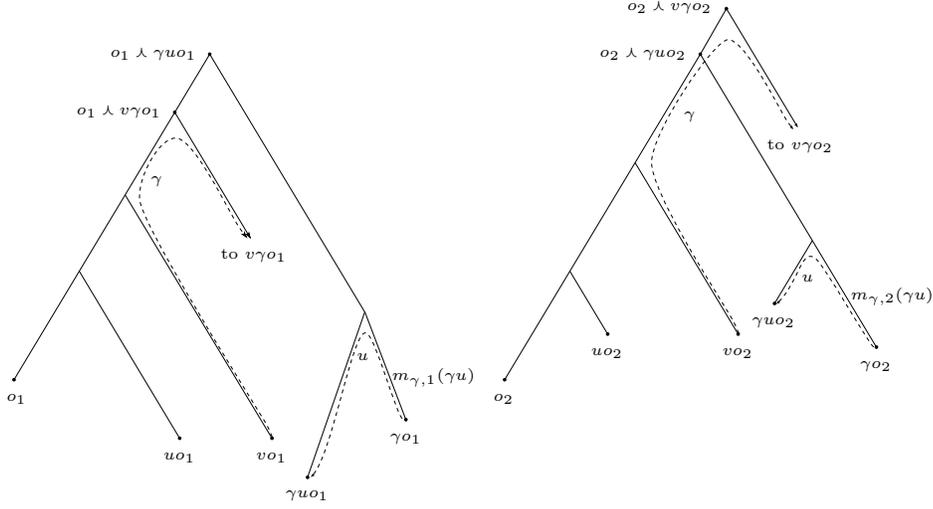
\begin{figure}[h!]
	\labellist \tiny \hair 5pt
	\pinlabel $o_1$ [t] at 5 121
	\pinlabel $uo_1$ [t] at 205 49
	\pinlabel $vo_1$ [t] at 318 49
	\pinlabel $\gamma o_1$ [t] at 483 72
	\pinlabel $\gamma uo_1$ [t] at 362 2
	\pinlabel $\textrm{to }v\gamma o_1$ [t] at 296 300
	\pinlabel $o_2$ [t] at 603 121
	\pinlabel $uo_2$ [t] at 730 177
	\pinlabel $vo_2$ [t] at 890 177
	\pinlabel $\gamma o_2$ [t] at 1060 161
	\pinlabel $\gamma uo_2$ [t] at 934 215
	\pinlabel $\textrm{to }v\gamma o_2$ [t] at 967 435
	\pinlabel $\gamma$ [l] at 154 367
	\pinlabel $u$ [t] at 430 171
	\pinlabel $\gamma$ [l] at 808 446
	\pinlabel $u$ [t] at 976 268
	\pinlabel $o_1\curlywedge \gamma uo_1$ [r] at 242 524
	\pinlabel $o_1\curlywedge v\gamma o_1$ [r] at 199 452
	\pinlabel $o_2\curlywedge \gamma uo_2$ [r] at 843 524
	\pinlabel $o_2\curlywedge v\gamma o_2$ [r] at 875 581
	\pinlabel $m_{\gamma,1}(\gamma u)$ at 517 125
	\pinlabel $m_{\gamma,2}(\gamma u)$ at 1080 225
	\endlabellist
	\centering \includegraphics[width=11.5cm]{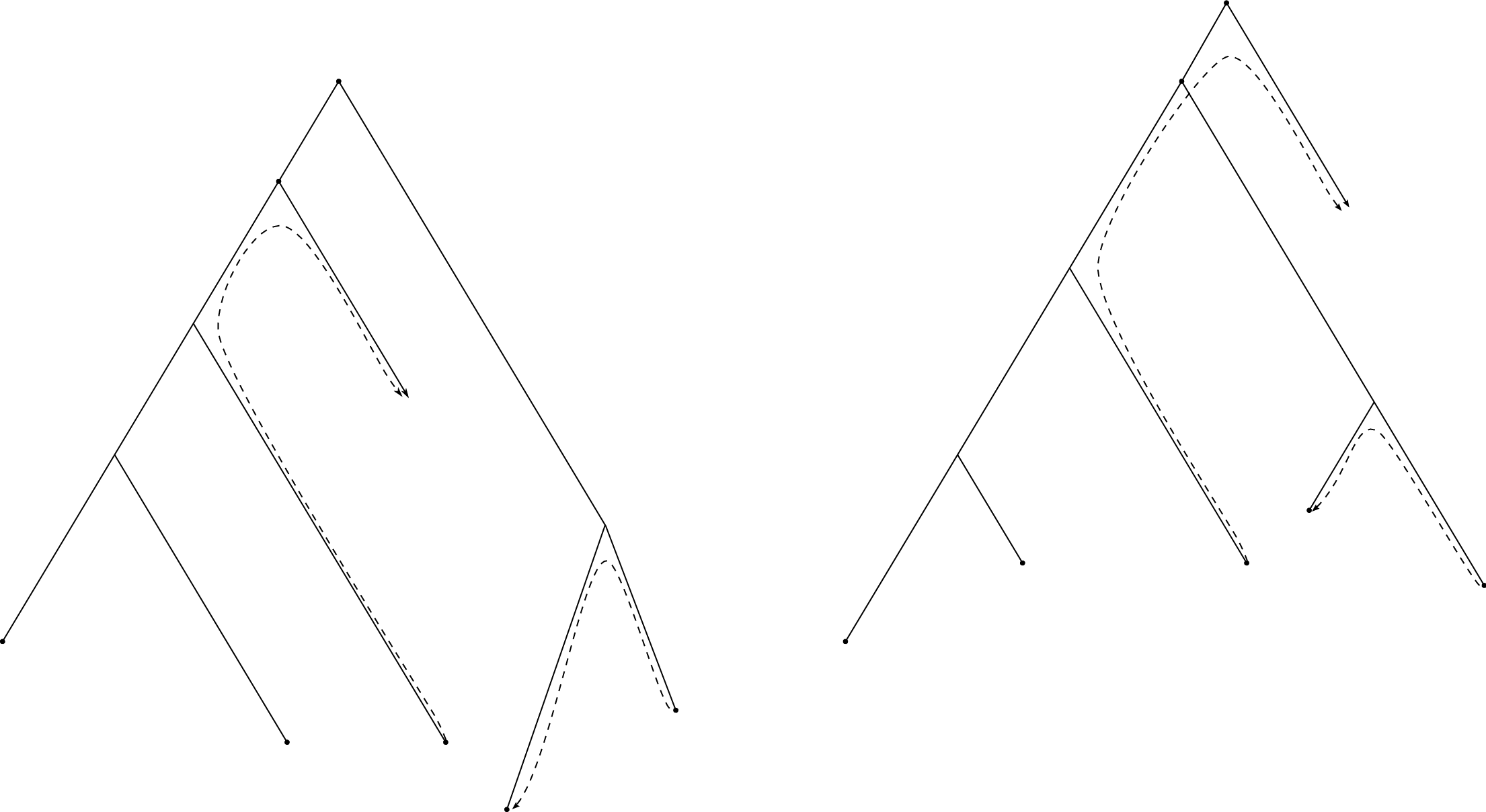}
	\caption[Case 2 of Theorem \ref{theorem:clf for lamplighter}]{The common ancestor $o_i \curlywedge \gamma u o_i$ lies in a different horocycle to $o_i \curlywedge v \gamma o_i$.} \label{fig:CLF lamp both trees}
\end{figure}

Suppose $m_i(\gamma) > m_i(u)+r$. 
If this is the case, we will have
	\begin{equation} \label{eq:anc1}
	o_i\curlywedge \gamma o_i = o_i \curlywedge \gamma u o_i
	\end{equation}
since the climb $m_{\gamma,i}(\gamma u)$ from $\gamma o_i$ to $\gamma u o_i \curlywedge \gamma o_i$ will be less than the fall $l_i(\gamma)$ from $o_i\curlywedge \gamma o_i$ to $\gamma o_i$. 
Indeed, the climb is given by $m_{\gamma,i}(\gamma u)$, which is equal to $m_i(u)$ by Lemma \ref{lem:m_i=m_h,i}.
By assumption $$m_i(u) < m_i(\gamma) -r < m_i(\gamma) -k.$$
Since $\h(\gamma o_i) = \pm k$, the fall $l_i(\gamma)$ from $o_i \curlywedge \gamma o_i$ to $\gamma o_i$ will be at least $m_i(\gamma) -k$. In summary, we get
		$$m_{\gamma,i}(\gamma u) \leq l_i(\gamma)$$
which implies \eqref{eq:anc1}.

Suppose $m_i(\gamma) > m_i(v) + r$.
In this case we claim 
	\begin{equation} \label{eq:anc2}
	o_i\curlywedge v \gamma o_i = v o_i \curlywedge v \gamma  o_i.
	\end{equation}
This follows from the fact that the climb $m_{v,i}(v\gamma)$ from $v o_i$ to $v o_i \curlywedge v\gamma o_i$ is longer than the fall $l_{i}(v)$ from $o_i \curlywedge v o_i$ to $v o_i$. Indeed, by Lemma \ref{lem:m_i=m_h,i} we have $m_{v,i}(v\gamma)=m_{i}(\gamma)$, which is greater than $m_i(v)+r$. But the fall $l_i(v)$ is at most the climb $m_i(v)+\h(v o_i) \leq m_i(v)+r$.

Thus, if we do not have \eqref{eqn:climbing limit} then we must have both \eqref{eq:anc1} and \eqref{eq:anc2}. However this implies
		$$o_i \curlywedge \gamma o_i = v o_i \curlywedge v\gamma o_i$$
which, by comparing their heights using Lemma \ref{lem:height of ancestor under action},  cannot occur when $r\neq 0$.
\end{proof}

\begin{proof}[Proof of Theorem \ref{thm:clf for lamplighter}]
With the structure of short conjugators understood from Proposition \ref{prop:short conjugators}, the linear upper bound when $r=0$ is clear since we may take $\gamma = (0,k)$ with $\modulus{k}\leq \modulus{u}+\modulus{v}=n$.

When $r\neq 0$, we use the conjguator $\gamma = (f,k)$ satisfying $0\leq k < \modulus{r}$ and \eqref{eqn:climbing limit}. Then Proposition \ref{prop:stein-taback} tells us that 
$$\modulus{\gamma} \leq 2m_1(\gamma)+2m_2(\gamma) \leq 2\max\{ m_1(u),m_1(v)\}+2\max\{ m_2(u),m_2(v)\}+4r.$$
As above, we may assume $r>0$. Since $\h(o_1\curlywedge u o_1) \leq -r$, we get $m_1(u) \geq r$, and similarly for $v$. Thus we will get
		\begin{eqnarray*}
			\modulus{\gamma} &\leq& 2 \max\{ m_1(u),m_1(v)\}+2\max\{ m_2(u),m_2(v)\} + 4\min\{m_1(u),m_1(v)\} \\
			&\leq & 3(m_1(u) + m_2(u) + m_1(v)+m_2(v))\\
			&\leq & 3 (\modulus{u}+\modulus{v}) = 3n.
		\end{eqnarray*}
This proves the Theorem.
\end{proof}

\section{A quadratic time algorithm}\label{sec:algorithm}

We finish by using the geometry developed in the preceding sections to describe an algorithm solving the conjugacy (search) problem in quadratic time.

Suppose $R$ is generated by a finite set $X$ as an abelian group. Then we can define a finite generating set for $\Gamma_2(R)$ as
		$$Y =\left\{ \left(\begin{array}{cc}t & x \\ 0 & 1\end{array}\right) : x \in X \right\}\mathrm{.}$$
We write the corresponding word-lengths as $\modulus{\cdot}_Y$ or $\modulus{\cdot}_S$, where $S$ is the generating set defined in Section \ref{sec:DL groups and gen lamplighters}. Note that when $R$ is finite we may take $X=R$ and then $Y=S$. We naturally get the inequality
		$$\modulus{u}_S \leq \modulus{u}_Y.$$

\begin{theorem}
	Suppose $R$ is finitely generated as an abelian group, with $X,Y$ as above. Then there is an algorithm which determines whether two elements $u,v$ in $\Gamma_2(R)$ are conjugate, and furthermore produces a conjugator, that runs in time $O(n^2)$, where $n=\modulus{u}_Y+\modulus{v}_Y$.
\end{theorem}

\begin{proof}
	We describe the steps of the algorithm below.
	
	\textbf{Step 1.}
	The input is given as words on the generating set $Y$.
	We may convert them into the form
			$u = (P,s),\ v=(Q,r)$
	in time linear in $n$. From this step we get a solution to the word problem in $\Gamma_2(R)$ which runs in linear time.
	
	\textbf{Step 2.}
	Check if $r=s$. If not, then stop and conclude that $u$ is not conjugate to $v$.
	If $r=s$, then continue to Step 3.
	
	\textbf{Step 3.}
	If $r=s\neq 0$ then proceed now to Step 4.
	Otherwise, by Proposition \ref{prop:short conjugators} part \eqref{part:unipotent}, there will be a conjugator $\gamma$ of the form
			$\gamma = (0,k)$
	with $\modulus{k} \leq \modulus{u}_S+\modulus{v}_S \leq n$. 
	Such $k$ must satisfy $P=t^kQ$.
	Calculate $k=v_0(P)-v_0(Q)$. 
	This can be done in linear time since we know $\modulus{v_0(P)},\modulus{v_0(Q)} \leq \modulus{u}_S+\modulus{v}_S \leq n$, so we find the minimal non-zero coefficient in $P$ or $Q$ of $t^i$ for $-n\leq i\leq  n$.
	If $\modulus{k}> n$ then we stop and conclude $u$ and $v$ are not conjugate.
	Otherwise we check whether $\gamma u=v \gamma$ using the linear time solution to the word problem.
	If it is then we stop the algorithm with the output that $u$ is conjugate to $v$ and $\gamma$ is the conjugator.
	If it is not then we stop and conclude that $u$ and $v$ are not conjugate.
	
	\textbf{Step 4.}
	We have $r=s\neq 0$. Without loss of generality we assume $r>0$, since otherwise we convert to $u^{-1}$ and $v^{-1}$.
	From Proposition \ref{prop:short conjugators}, we know that if $u,v$ are conjugate then there is a conjugator $\gamma$ of the form
			$\gamma = (f,k)$
	where $0\leq k < r$. Since the climbing of $\gamma$ is limited by \eqref{eqn:climbing limit}, we know we can write $f$ in the form
			$$f = \sum\limits_{i=n_1}^{n_2} a_it^i \ \textrm{ with } \  \modulus{n_i}= \max\{m_i(u),m_i(v)\} + r \leq 2n.$$
	Rearranging equation \eqref{eq:lamplighter functions}, we get $t^rf - f = t^kQ - P$. 
	For each $k=0,1,2,\ldots,r-1$, do the following:
	let 
			$$t^kQ-P = \sum\limits_{i=-\infty}^\infty b_i t^i.$$
	By comparing the coefficients of $t^{n_2+1},\ldots,t^{n_2+r}$ in $t^rf - f = t^kQ - P$, we may set $a_{i} = b_{i+r}$ for $i=n_2-r+1,\ldots, n_2$.
	Now we can compare the coefficients of $t^{n_2-r+1},\ldots,t^{n_2}$, and thus set $a_i = a_{i+r} + b_i$ for $i = n_2-2r+1,\ldots,n_2-r$. 
	We repeat this process until all coefficients $a_i$ have been assigned, and we check whether $a_i=b_i$ for $i=n_1,\ldots, n_2-\alpha r$, where $\alpha$ is the maximal integer such that $n_1 \leq n_2-\alpha r$.
	
	Since $n_2-n_1 \leq 4n$, the process of assigning elements of $R$ to each $a_i$ will run within $O(n)$ time.
	If $a_i=b_i$ for $i=n_1,\ldots, a_{n_2-\alpha r}$ then we stop and $\gamma$ is a conjugator for $u,v$.
	If for each $k$, there is some $i \in \{n_1,\ldots,n_2-\alpha r\}$ for which $a_i \neq b_i$, then we stop and conclude that $u,v$ are not conjugate. We have to do this at most $r$ times, so in total the fourth step runs within quadratic time.
\end{proof}

\textbf{Acknowledgements:}
Thanks are due to C.{} Abbott, Y. Antol\'in,  C.{} Dru\c{t}u, \mbox{T.{} Riley} and an anonymous referee for helpful discussions and comments.

\bibliographystyle{amsplain}
\bibliography{bibliography}

\end{document}